\documentclass[12pt]{article}
\usepackage{amsmath,amsfonts,amssymb,amsthm}

\def\noi{\noindent}

\def\pf{\noi{\bf Proof.\ \,}}
\def\eop{{$\square$}}

\def\labtt#1{\label {#1}} 
\def\labttr#1{\label {#1} \rm }

\def\vep{\varepsilon}

\def\FF{{\mathbb F}}
\def\QQ{{\mathbb Q}}

\def\ZZ{{\mathbb Z}}

\def\la{\langle}
\def\ra{\rangle}
\def\<{\langle}
\def\>{\rangle}

\def\half{{1 \over 2}}

\def\eighth{{1 \over 8}}

\def\explus#1{2^{1+2#1}_+}

\def\dg#1{{\cal D}({#1})}  

\def\gd#1{G_{2^{#1}}} 
\def\rd#1{R_{2^{#1}}}

\begin{document}

\newtheorem{thm}{Theorem}[section]
\newtheorem{prop}[thm]{Proposition}
\newtheorem{lem}[thm]{Lemma}
\newtheorem{rem}[thm]{Remark}
\newtheorem{coro}[thm]{Corollary}
\newtheorem{conj}[thm]{Conjecture}
\newtheorem{de}[thm]{Definition}
\newtheorem{hyp}[thm]{Hypothesis}

\newtheorem{nota}[thm]{Notation}
\newtheorem{ex}[thm]{Example}
\newtheorem{proc}[thm]{Procedure}

\centerline{9 November, 2005}
\begin{center}
{\Large   Corrections and additions to `` Pieces of $2^d$: existence and uniqueness for
Barnes-Wall and Ypsilanti lattices. '' }

\vspace{4mm}
Robert L.~Griess Jr.
\\[0pt]
Department of Mathematics\\[0pt] University of Michigan\\[0pt]
Ann Arbor, MI 48109-1043  \\[0pt]
\vskip 1cm 

\end{center}

\begin{abstract}  Mainly, we correct the uniqueness result by adding a projection requirement to condition X 
and  give a better proof for the equivalence of commutator density, 2/4-generation and 3/4-generation.
\end{abstract} 

\tableofcontents

\section{The main results for for Barnes-Wall lattices in \cite{bwy}}  

{\it
The single change for Section 2 of \cite{bwy} is to add condition (f) to Condition $X(2^d)$, 2.3.  The corrected hypothesis is stated below.  
}

\begin{de}\labttr{condx} {\bf Condition $X(2^d)$: }
This is defined for integers $d\ge 2$.  Let $s \in \{0,1\}$ be the
remainder of $d+1$ modulo 2. 

We say that the quadruple $(L,L_1,L_2,t)$    
is a {\it an X-quadruple} if it 
satisfies {\it condition $X(2^d)$} (or, more simply, 
{\it condition X}), listed below:

(a)  $L$ is a rank $2^d$ even integral lattice containing 
$L_1\perp L_2$, 
the  orthogonal direct sum of 
sublattices $L_1\cong L_2$ of rank
$2^{d-1}$;

(b) When $d=2$, $L\cong L_{D_4}\cong BW_4$ 
and $L_1\cong
L_2\cong L_{A_1^2}$; when $d\ge 3$, 
$2^{-{s\over 2}}L_1$ and $2^{-{s\over 2}} L_2$ 
are initial entries of quadruples which satisfy
condition  $X(2^{d-1})$.  

(c) $\mu (L) = 2^{\lfloor {d\over 2} \rfloor }$.

(d)   $\dg L
\cong 2^{2^{d-1}}, 1$ as $d$ is even, odd, respectively.  

(e) There is an isometry $t$ of order 2 on $L$ 
which interchanges $L_1$ and $L_2$ and 
satisfies $[L,t]\le L_1\perp
L_2$, i.e., acts trivially on $L/ [L_1 \perp L_2 ]$.  

(f)  The projection of $L$ to each $V_i$ are rank $2^{d-1}$ BRW-sublattices, i.e., are stable under 
a natural $BRW^0(2^{d-1},+)$-subgroup of $Aut(L)$ which centralizes $t$ and stabilizes each $L_i$. 
\end{de}

\section{Revision of equivalence of commutator density, 2/4- and 3/4-generation}

{\it 
This section should replace Section 5.1 of \cite{bwy}(items 5.16 to 5.22), which shows equivalence of commutator density, 2/4-generation and 3/4-generation.  We recall the definitions here.   Let $M$ be a module for a dihedral group $D$ of order 8 so that the central involution of $D$ acts on $M$ as $-1$.  Let $f \in D$ be an element of order 4.  Commutator density is the property that $[M,D]=M(f-1)$.  The 3/4 generation property is the condition that $M$ is generated by fixed points of any three noncentral involutions of $D$.  The 2/4 generation property is the condition that $M$ be generated by the fixed points of a generating pair of involutions in $D$.

To fit with \cite{bwy}, change numbers in the revision below to begin with ``5.16".  } 

\bigskip 

In this section, ``lattice'' means just a free abelian group since the bilinear form is irrelevant to the arguments.  

\begin{nota}\labttr{nota0} Let $L$ be a lattice  with involution $s$ acting on $L$.  
Let $L^\vep (s):=\{ v \in L | vs=\vep v\}$ be the $\vep$-eigenlattice for the involution $s$ and let $Tel(s):=L^+(s)\oplus L^-(s)$ be the total eigenlattice for $s$.  
\end{nota}

\begin{nota}\labttr{revnota1}
Let 
$t, u$ be involutions which generate the dihedral group $D$ of order 8 and let $f$ be an element of order 4 in $D$.  
\end{nota}

\begin{nota}\labttr{nota1}
Let  $L$ be  a free abelian group of rank $2n$ which admits an action of $D$ in which the central involution $[t,u]$ acts as $-1$ 
and let $L_1:=L^-(t), L_2:=L^+(t)$.   

We extend the action of $D$ to the ambient vector space $\QQ \otimes L$.  
For integers $\ell \le m$, let $Q(\ell, m):=2^\ell Tel(t) / 2^m Tel(t)$.  

Define integers $d, e$ to be the number of Jordan blocks of size 2, 1 respectively, for the action of $u$ on the elementary abelian group 
$L/Tel(t)$.  
\end{nota}

\begin{lem}\labtt{ltel}
$|L:Tel(t)|=2^{2d+e}$ and 
$|L:2Tel(t)|=2^{2n+2d+e}$.  
\end{lem}

\begin{nota}\labttr{revnota2} We have a chain of $D$-invariant abelian groups  
$2Tel(t) \le L \le \half Tel(t)$.  For $g\in D$, we denote by $A(g)$ and  $B(g)$ the commutator modules $Q(-1,1)(g-1) = \half Tel(t)(g-1) + 2Tel(t)/2Tel(t)$ and  $L(g-1)  + 2Tel(t)/2Tel(t)$, respectively. 
\end{nota}

\begin{lem} \labtt{intaq0} 
$A(u)\cap Q(0,1) = B(u)\cap Q(0,1)=Q(0,1)(u-1)$.  
\end{lem} 
\pf  Clearly, $A(u)\cap Q(0,1) \ge  B(u)\cap Q(0,1)\ge Q(0,1)(u-1)$.  
Now to prove the opposite containment.  Since $A(u)\cap Q(0,1)$ consists of elements inverted by $u$, hence fixed by $u$, it is contained in the subgroup  $Q(0,1)(u-1)$ of the free $\FF_2\la u\ra$ module $Q(0,1)$.  \eop

\begin{lem}\labtt{commsum}   $A(t)\cap A(u)=0$.  
\end{lem}
\pf Since $Tel(t)(t-1)=2L_1$, 
the image of $A(t)$ is just $L_1+2Tel(t)/2Tel(t)\le Q(0,1)$.   
Also, $A(u) \cap Q(0,1)$ is exactly the image of the diagonal sublattice $\{ (x, xu)|x\in L_1\}$ 
of $L_1\oplus L_2$ in $Q(0,1)$.  The result follows.  
\eop

\begin{lem}\labtt {liftfixedpoints}  
A coset of $Tel(t)$  fixed by $u$ contains an element fixed by $u$. 
\end{lem} 
\pf  Let $x+Tel(t)$ be such a coset.  
Since $Tel(t)$ is a free $\ZZ \la u \ra$-module, every element of $Tel(t)$ negated by $u$ is a commutator.  
Therefore, there exists $v \in Tel(t)$ so that $x(u-1)=v(u-1)$.  
Then $x-v$ is fixed by $u$ and is in $x+Tel(t)$. 
\eop

\begin{lem}\labtt{orderbu} $|B(u)|=2^{n+d}$.
\end{lem}
\pf
The right side is the product of $|Q(0,1)|^{\half}=2^n$ with $2^d=|(L/Tel(t))(u-1)|$. 
To evaluate the left side, 
use 
\ref{intaq0}, 
\ref{liftfixedpoints}.
\eop 

\begin{lem} \labtt{kerendo} The kernel of the endomorphism induced by $t-1$ on $L/2Tel(t)$ 
is just $Q(0,1)$.
\end{lem}
\pf
If the kernel were larger, there would be 
$x\in L\setminus Tel(t)$ 
so that
$x(t-1)\in 2Tel(t)$.  
Then there would be a unique $y \in L^-(t)$ 
so that $2y=x(t-1)$, 
whence 
$yt=-y$ and 
$(x+y)t=x+y$ and so $x+y\in L^+(t)$, 
a contradiction to $x \not \in Tel(t)$.  
\eop

\begin{coro}\labtt{orderbt} $|B(t)|=2^{2d+e}$.
\end{coro}
\pf
Since $Q(0,1)$ has index $2^{2d+e}$ in $L/2Tel(t)$, the result follows from \ref{kerendo}.   
\eop

\begin{lem}\labtt{orderbtbu} 
$B(t) \cap B(u)=0$ and $|B(t)+B(u)|=2^{n+3d+e}$.  
\end{lem} 
\pf  For $B(t)\cap B(u)=0$, use \ref{commsum}.  
The second statement follows from the formula $|B(t)+B(u)|=|B(t)||B(u)|/|B(t)\cap B(u)|=|B(t)||B(u)|$ and \ref{orderbu}, \ref{orderbt}.  
\eop 

\begin{lem}\labtt{tuf}$L(t-1)+L(u-1)=[L,D]\ge L(f-1) \ge 2L$.  
\end{lem}
\pf  First, $L(t-1)+L(u-1)=[L,D]$ holds 
because $t$ and $u$ generate $D$.  
The containment $[L,D]\ge L(f-1)$ is obvious.  
Since $[L,D]$ contains $L(f-1)$ and $L(f-1)^2=2L$, the final containment holds. 
\eop

\begin{lem}\labtt{cd=3/4}  We have 
$|B(t)+B(u):(L(f-1)/2Tel(t))|=2^{n+3d+e-(2d+e+n)}=2^d$.  
Therefore, 
$L(f-1)=L(t-1)+L(u-1)$ if and only if $d=0$.  In other language, commutator density is equivalent to 2/4 generation.  
\end{lem} 
\pf 
Observe that if $M$ is any $f$-invariant subgroup of $L$, then $M(f-1)^2=2M$ and that for any integer $j$, 
(*) 
$|M(f-1)^j:M(f-1)^{j+1}|=2^n$.  Since $Tel(t)(f-1)\le Tel(t)$, we have  $2Tel(t)=Tel(t)(f-1)^2 \le L(f-1)\le [L,D]$.  Both $L(f-1)$ and $L(t-1)+L(u-1)=[L,D]$ contain $2Tel(t)$, whence a basic isomorphism theorem implies that $|L(t-1)+L(u-1):L(f-1)|=|B(t)+B(u):(L(f-1)/2Tel(t))|$.  
The statements follow from \ref{ltel}, \ref{orderbtbu} and (*).  \eop

\begin{lem}\labtt{3/4=2/4} 
The properties 2/4-generation and 3/4-generation are equivalent.  
\end{lem} 
\pf
Obviously, 2/4-generation implies 3/4-generation.  Conversely, assume that  $L=L^+(t)+L^-(t) +L^+(u)$.  Using \ref{tuf}, 
we have $L^+(t)+L^+(u)\ge L(t+1)+L(u+1)= [L,D]\ge 2L$.  Since $L^+(t)+2L=L^-(t)+2L=L^+(t)+L^-(t)$,  $L=L^+(t)+L^-(t) +L^+(u)=L^+(t)+L^+(u)+2L=L^+(t)+L^+(u)$, whence 2/4-generation.  
\eop  

Lemmas \ref{cd=3/4} and \ref{3/4=2/4} imply the required equivalence of commutator density, 2/4-generation and 3/4-generation.  

\begin{rem}\rm Note that the present version avoids bilinear forms, so does not use the language of determinants.   If a $D$-invariant bilinear form is present, we may replace statements about finite indices of sublattices with ones about determinants.  Note that, for a noncentral involution $s$ of $D$, each $L^\pm (s)$ has the same determinant.  
\end{rem}

\section{Revision of proof of the uniqueness theorem}

{\it The following is a replacement for Section 9 in \cite{bwy}  
``Proof of Uniqueness''.   The statement of Theorem 9.2 is unchanged except for correcting a quotient group in (ii)(a), including a definition in  (ii)(c) and converting its final sentence to a remark.  The proof is revised, using $X(2^d)(f)$.   

Item numbers should be changed to begin with ``9.1''.  Numbered references internal to the proof refer to \cite{bwy}. }

\begin{nota}\labttr {frakx} 
Given $d\ge 3$ and  $L_1, L_2$, we let 
${\frak X} := {\frak  X}(L_1,L_2)$ be the set of all
$X$-quadruples of the form $(L,L_1,L_2,t)$;  see 2.3.  
\end{nota}

\begin{thm} \labtt{bwtclassification} 
We use the notation in 2.3, 6.1, 6.2, 7.2 and 9.1.  
Suppose that $d \ge 3$ and $(L_1, L_2)$ is an orthogonal pair of
lattices, so that each 
$L_i$ is BW-type of rank $2^{d-1}$.  

(i)  ${\frak X}$ is an 
orbit under the natural action of $F_1 \times F_2$, where
$F_i:=Stab_{Aut(L_i)}(L_i[1-r])$ (see 6.2; by
induction,
$F_i
\cong \gd {d-1}$).  Define $Q_i:=C_{F_i}(L_i/L_i[1])$.  

The elements of $\frak X$ are in correspondence with each of the following 
sets.

(a) $F_1/Q_1$;  

(b) $F_2/Q_2$;

(c) Pairs of involutions $\{s, -s\}$ in the orthogonal group on $V$ which
interchange
$L_1$ and $L_2$.

(d) Dihedral groups of order $8$ which are generated by the SSD
involutions associated to $L_1, L_2$ and involutions as in (c).

(ii) (a) The subgroup $G_L^0$ of $F_1 \times F_2$  which stabilizes
$L$ has structure 
$Q_1 \times Q_2 \le G_L^0$ and $G_L^0/(Q_1\times Q_2)$ is the diagonal subgroup of 
$F_1/Q_1 \times F_2/Q_2$ with respect to the isomorphism induced by
$s$, an involution as in (i.c).

(b) 
The subgroup $G_L$ of $Aut(L_1 \perp L_2)
\cong Aut(L_i)
\wr 2$ which stabilizes $L$ is
$G_L^0\la  s \ra$.  
We have $G_L \cong [\explus {(d-1)} \times \explus
{(d-1)}].[\Omega^+(2(d-1),2)
\times 2]$.  

(c) The subgroup of $G_L$ which acts trivially on $L/L[1]$ is $R:=\la
Q_{12}, s, t_i \ra$, where $t_i$ is the SSD involution associated to
$L_i$ and $Q_{12}:=\{xx^s | x \in Q_1\}$.  The quotient
$G_L/R
\cong 2^{2d-2}{:}\Omega^+(2d-2,2)$
is a maximal parabolic subgroup of $Out^0(\explus {d})
\cong
\Omega^+ (2d,2)$. (See Appendix A.0).  
\end{thm}

\begin{rem} \rm The extension in 
(c) is split, despite
$\gd e$ being nonsplit over $\rd e$ for $e \ge 4$.  See Appendix A2.  
\end{rem}

\pf  (i) 
We prove the classification by induction.  For $d=2$,
$Aut(L_{D_4})\cong  2^{1+4}_+[Sym_3 \wr 2]$ and for $d=3$, 
$Aut(L_{E_8} )\cong
W_{E_8}$.  
When $d=4$, the main theorem follows from the arguments of \cite{POE}.  

For the rest of the proof, we assume that $d \ge 4$.  
By induction, a lattice satisfying the $X(2^{d-1})$ condition is uniquely
determined up to isometry.  
This applies to the lattices $L_1,  L_2$.

Let $(L, L_1, L_2, t) \in {\frak X}$.  Then $det(L)$ and $|L:L_1\perp L_2|$ are determined.  Define $G:=Aut(L)$.

Let $p_i$ be the orthogonal projection of $L$ to $V_i:=\QQ\otimes L_i$, for $i=1,2$.  
Since $L$ is integral, 
$L^{p_i}$, the projection of $L$ to $V_i$ is 
contained in $L_i^*$.  
When $d$ is odd, $L^{p_i}$ must be $L_i^*=L_i[-1]$, by determinant  considerations.  

Assume now that $d$ is even.  
There is a subgroup $H$ of $C_G(t)$, $H\cong BRW^0(2^{d-1},+)$, so that $H$ acts faithfully on both $V_i$ and stabilizes  $L$, $L_1$ and $L_2$.  
By 7.14 and $X(2^d)(f)$, the projection is $L_i[-1]$,
where the twist is with respect to a lower fourvolution in $H$.    
Since the groups $Q_i$ act trivially on $L_i[-1]/L_i$, for $i=1,2$, the group $\la Q_1, Q_2, t \ra$ stabilizes $L$ and in fact acts trivially on $L/L_1\perp L_2$.  
We  define a group $R:=\la Z(Q_1), Z(Q_2), Q_{12}, t\ra \cong \explus d$, where $Q_{12}:=\{ xx^t|x \in Q_1\}$.

Let $D$ be the dihedral group which is generated by $t$ and either $t_i$, the involution generating $Z(Q_i)$ (or what is the same, the RSSD involution associated to $L_i$).  
It follows that $L$ is determined by $L_1, L_2, D$ 
in the sense that 
$L$ lies between $L_1\perp L_2$ and 
$L_1[-1]\perp L_2[-1]$ and 
$L/[L_1\perp L_2]$ is the fixed point submodule 
for the action of $D$ (equivalently, of $t$) on $[L_1[-1]\perp L_2[-1]]/[L_1\perp L_2]$.  Recall from earlier in the proof that $D$ 
does determine the $L_i[-1]$ by use of X(f).

Now, to what extent do $L$ and $L_1\perp L_2$ determine $D$?  The answer is: up to conjugacy in $Aut(L_1\perp L_2)\cong Aut(L_1)\wr 2$ (note that $d\ge 4$ here).  
Our group $D$ is generated by the center of the natural index 2 subgroup
of $Aut(L_1 \perp L_2)$ and a wreathing involution.  
In general, 
wreathing involutions in a wreath product of groups $K \wr 2$
form an orbit under the action of either direct factor isomorphic to $K$
in the base group of the wreath product. This proves correspondence of $\frak X$ with
(c) and (d).  The stabilizer subgroup is diagonal in the base group $K
\times K$, and either direct factor represents all cosets of the stabilizer
(whence the equivalence of $\frak X$  with (a) and (b))

It follows that, up to isometry preserving $L_1 \perp L_2$,
$D$, hence $L$, is determined by the pair of indecomposable lattices
$L_1$ and
$L_2$.

Proof of statement (ii) is easy.  The statement about
parabolic subgroups is proven with a standard result from  
the theory of
Chevalley groups, e.g.
\cite {Car}.  Independently of that theory, the maximality could be proved
directly by showing that there is no system of imprimitivity on the set
of isotropic points.    This is an exercise with Witt's theorem. 
\eop

\section{Other Revisions}

In the Abstract, p.147, line 10, insert ``which'' after ``$(log_2(mass(n)))$''.

End of first paragraph Theorem 2.7.  Change 
``$(\eighth + o(1))d \, 2^d$.'' to 
``$(\eighth + o(1))d \, 2^d$).''

In Appendix 2 of \cite{bwy}, the statement of Lemma 15.9 is incorrect (and the lemma is not used).  A correct analysis of involutions in the BRW-group is in \cite{ibw1}.

\end{document}